\documentclass[11pt, a4paper, reqno]{amsart}
\usepackage[utf8]{inputenc}
\usepackage{amsmath,amssymb}
\usepackage{graphicx}
\usepackage{indentfirst,csquotes}
\usepackage[colorlinks=true, linkcolor=blue, citecolor=blue, urlcolor=blue]{hyperref}
\usepackage{pdflscape} 
\usepackage{hhline} 
\usepackage{array} 
\usepackage{calrsfs}
\DeclareMathAlphabet{\pazocal}{OMS}{zplm}{m}{n}
\usepackage{latexsym,amsmath,amsfonts,amscd,amssymb, mathrsfs, amsthm}
\usepackage[english]{babel}
\usepackage{appendix}
\usepackage{hyperref}
\usepackage{mathtools}
\usepackage{booktabs}
\usepackage{natbib}

\advance\oddsidemargin by -0.5cm
\advance\evensidemargin by -0.5cm
\advance\topmargin by -1.5cm
\def\mxl{\left[ \begin{array}}  \def\mxr{\end{array} \right]}
\def\detl{\left| \begin{array}}  \def\detr{\end{array} \right|}
\def\bmx{\begin{bmatrix}} \def\emx{\end{bmatrix}}
\def\mx#1{\begin{bmatrix}#1\end{bmatrix}}

\def\bc{\begin{center}}          \def\ec{\end{center}}
\def\ben{\begin{enumerate}}      \def\een{\end{enumerate}}
\def\beq{\begin{equation}}       \def\eequ{\end{equation}}
\def\beqn{\begin{eqnarray*}} \def\eeqn{\end{eqnarray*}} 
\def\eqs{\vspace{-10pt}\begin{equation}} \def\eqe{\vspace{-5pt}\end{equation}}
\def\bal{\begin{align}}           \def\eal{\end{align}}
\def\ben{\begin{enumerate}}       \def\een{\end{enumerate}}
\def\bit{\begin{itemize}}           \def\eit{\end{itemize}}
\def\btabb{\begin{tabbing}}         \def\etabb{\end{tabbing}}
\def\btab{\begin{tabular}}          \def\etab{\end{tabular}}
\def\barr{\begin{array}}           \def\earr{\end{array}}
\def\earrb{\end{array} \right\}}
\def\xx{\times}

 \def\gl{\lambda}

\def\it{\emph}  
\def\bf{\textbf}

\def\0{\bf{0}}    
\font\Bbb=msbm10 
\def\A{\ifmmode{\mbox{\Bbb A}}\else{{\Bbb A}}\fi}
\def\B{\ifmmode{\mbox{\Bbb B}}\else{{\Bbb B}}\fi}
\def\R{\ifmmode{\mbox{\Bbb R}}\else{{\Bbb R}}\fi}
\def\C{\ifmmode{\mbox{\Bbb C}}\else{{\Bbb C}}\fi}

\def\G{\ifmmode{\mbox{\Bbb G}}\else{{\Bbb G}}\fi}
\def\N{\ifmmode{\mbox{\Bbb N}}\else{{\Bbb N}}\fi}
\def\Z{\ifmmode{\mbox{\Bbb Z}}\else{{\Bbb Z}}\fi}
\def\K{\ifmmode{\mbox{\Bbb K}}\else{{\Bbb K}}\fi}
\def\L{\ifmmode{\mbox{\Bbb L}}\else{{\Bbb L}}\fi}
\def\D{\ifmmode{\mbox{\Bbb D}}\else{{\Bbb D}}\fi}
\def\T{\ifmmode{\mbox{\Bbb T}}\else{{\Bbb T}}\fi}
\def\I{\ifmmode{\mbox{\Bbb I}}\else{{\Bbb I}}\fi}
\def\E{\ifmmode{\mbox{\Bbb E}}\else{{\Bbb E}}\fi}
\def\K{\ifmmode{\mbox{\Bbb K}}\else{{\Bbb K}}\fi}
\def\Q{\ifmmode{\mbox{\Bbb Q}}\else{{\Bbb Q}}\fi}
\def\P{\ifmmode{\mbox{\Bbb P}}\else{{\Bbb P}}\fi}
\def\cG{\ifmmode{\mathcal{G}}\else{{$\mathcal{G}$}}\fi}
\def\cA{\ifmmode{\mathcal{A}}\else{{$\mathcal{A}$}}\fi}

\def\op{ \ifmmode{\oplus} \else{\leavevmode\hbox{$\oplus$}\fi } }
\def\ot{\ifmmode{\otimes}\else{$\!\!\otimes$}\fi}
\def\od{ \ifmmode{\odot} \else{\leavevmode\hbox{$\odot$}\fi } }
 
\def\adots{\mathinner{\raise 1pt\hbox{.}\mkern1mu\raise4 pt\hbox{.}\mkern2mu
  \mkern1mu\raise7 pt\vbox{\kern7 pt\hbox{.}}\mkern2mu}}

\def\bdefn{\begin{defn}} \def\edefn{\end{defn}}
\def\bex{\begin{ex}} \def\eex{\end{ex}}
\def\beg{\begin{eg}} \def\eeg{\end{eg}}
\def\blem{\begin{lem}} \def\elem{\end{lem}}
\def\bth{\begin{thm}} \def\eth{\end{thm}}
\def\bprop{\begin{prop}} \def\eprop{\end{prop}}
\def\bcor{\begin{cor}} \def\ecor{\end{cor}}

 \def\dag{\dagger}
 
\def\FF{\mathbb{F}}
 \newtheorem{thm}{Theorem}[section]
 \newtheorem{cor}[thm]{Corollary}
 \newtheorem{lem}[thm]{Lemma}
 \newtheorem{prop}[thm]{Proposition}
 \theoremstyle{definition}
 \newtheorem{defn}[thm]{Definition}
 \theoremstyle{remark}
 
 \newtheorem*{ex}{Example}
 \numberwithin{equation}{section}

\def\ball{\begin{align}}  \def\eall{\end{align}}

\def\bbin{\left( \begin{array}{c}} \def\ebin{\end{array} \right)} 

\begin{document}
\title[On generalized inverses of graph matrices]{On generalized inverses of matrices associated with certain graph classes}
\author[Claúdia M. Araújo]{Cláudia Mendes Araújo}
\author[Faustino A. Maciala]{Faustino A. Maciala}
\author[Pedro Patrício]{Pedro Patrício}
\date{\today}
\address{Claúdia Mendes Araújo, CMAT -- Centro de Matemática and Departamento de Matemática Universidade do Minho 4710-057 Braga Portugal}
\email{clmendes@math.uminho.pt}
\address{Faustino A. Maciala, CMAT -- Centro de Matemática, Universidade do Minho, 4710-057 Braga, Portugal; 
Departamento de Ciências da Natureza e Ciências Exatas do Instituto Superior de Ciências da Educação de Cabinda -- ISCED-Cabinda, Angola.}
\email{fausmacialamath@hotmail.com}
\address{Pedro Patrício, CMAT -- Centro de Matemática and Departamento de Matemática Universidade do Minho 4710-057 Braga Portugal}
\email{pedro@math.uminho.pt}
\maketitle
\let\thefootnote\relax
\footnotetext{MSC2020: Primary 15A09, Secondary 05C20, 05C50.} 

\begin{abstract}
We investigate generalized inverses of matrices associated with two classes of digraphs: double star digraphs and D-linked stars digraphs. For double star digraphs, we determine the Drazin index and derive explicit formulas for the Drazin inverse. We also provide necessary and sufficient conditions for the existence of the Moore--Penrose inverse and give its explicit expression whenever it exists. For $D$-linked stars digraphs, we characterize when the group inverse exists and obtain its explicit form. In the singular case where $BC = 0$, we express the Drazin index of the matrix in terms of the Drazin index of the base digraph matrix. Additionally, we establish necessary and sufficient conditions for Moore--Penrose invertibility and derive explicit formulas in that case. Our results reveal a clear connection between the algebraic structure of generalized inverses and the combinatorial properties of these graph classes, providing a unified framework for group, Drazin, and Moore--Penrose invertibility.
\end{abstract} 

\bigskip
\noindent 
\section{\textbf{Introduction}}
\label{sec:introduction}
\hspace{0.5cm} A central problem in spectral graph theory and matrix analysis is to determine when graph-associated matrices are invertible and, whenever possible, to derive explicit formulas for their inverses. For example, an explicit inverse formula for the distance matrix of a wheel graph was obtained in \cite{Balaji}, while the invertibility of trees and bipartite graphs was investigated in \cite{Godsil} and \cite{Bapat, McLeman}, respectively. The latter results were later extended to the framework of Drazin invertibility in \cite{Catral2}. Both group and Drazin inverses have been extensively studied for matrices associated with special classes of graphs. In particular, the group inverse has been characterized for path graphs, tree graphs, and bipartite digraphs in \cite{Catral-grp}. Recent progress has also been made on generalized inverses of graph matrices, for instance in \cite{Abiad}, where the group inverse of $M$-matrices associated with distance-biregular graphs was characterized. Such developments further illustrate the deepening interplay between algebraic properties of matrices and the combinatorial structure of graphs -- a connection that also underpins the present work. For a matrix-oriented approach to graph theory, we follow the framework established in \cite{bookBapat}. 
Beyond their intrinsic algebraic relevance, generalized inverses have a broad range of applications. The Drazin inverse plays a crucial role in the analysis of Markov chains, singular differential equations, and network dynamics, whereas the Moore--Penrose inverse arises naturally in optimization, data science, and control theory. In graph-theoretic contexts, these inverses capture structural properties of graphs and digraphs, linking combinatorial features to algebraic and spectral invariants. In particular, the Drazin index quantifies the level of singularity of a matrix and serves as a measure of its deviation from invertibility.

The present work advances the study of generalized inverses of graph matrices by addressing in detail the Drazin and Moore--Penrose invertibility of matrices associated with two families of digraphs -- the double star digraphs and the $D$-linked stars digraphs. For double star digraphs, we determine the Drazin index and derive explicit formulas for the Drazin inverse, as well as for the Moore--Penrose inverse whenever it exists, providing necessary and sufficient conditions for Moore--Penrose invertibility. For D-linked stars digraphs, we characterize the existence of the group inverse and obtain its explicit expression; in the singular case where \(BC = 0\), we determine the Drazin index of \(M\) in terms of that of the underlying digraph matrix \(A\). We also provide explicit formulas for the Moore--Penrose inverse of \(M\), together with necessary and sufficient conditions for its existence. These results strengthen the link between the algebraic properties of generalized inverses and the combinatorial structure of graphs, integrating group, Drazin, and Moore--Penrose inverses within a unified analytical framework. The explicit block formulas further enhance applicability to computational and algorithmic settings, including network analysis, structured optimization, and singular dynamical systems.

The paper is organized as follows. Section \ref{sec:preliminaries} establishes the algebraic framework underlying our results, recalling the main definitions and properties of the Drazin and Moore--Penrose inverses. Special attention is given to the Drazin inverse, whose existence and index are characterized in terms of the minimal polynomial of the matrix. Section \ref{SecDoubleStar} focuses on matrices associated with double star digraphs, determining their Drazin index in all non-invertible cases and deriving explicit inversion formulas. Section \ref{SecLinked} extends the analysis to $D$-linked stars digraphs, providing conditions for group, Drazin, and Moore--Penrose invertibility. The paper concludes with remarks and open problems that suggest further directions in the study of generalized inverses of graph matrices.

\section{\textbf{Preliminaries}}
\label{sec:preliminaries}
\hspace{0,5cm} For completeness, we briefly recall the main definitions and properties of generalized inverses -- particularly the Drazin and Moore--Penrose inverses -- which provide the algebraic framework for our analysis.

Throughout this paper, we work over a field $\FF$, possibly endowed with an involution $(\bar{\cdot})$. In this setting, the involution in $\FF$ induces an involution $*$ on the space of finite matrices over $\FF$, defined by $A^*=[a_{ij}]^*=[\overline{a_{ji}}]$. Also, $\overline{A}=\overline{[a_{ij}]}=[\overline{a_{ij}}].$  

The Drazin inverse of an $n \times n$ matrix $A$, denoted by $A^D$, is the unique matrix satisfying
\[
A^{k+1}X = A^k, \quad XAX = X, \quad AX=XA,
\]
for some nonnegative integer $k$. The smallest such $k$ is called the \emph{Drazin index} of $A$, denoted by ${\rm ind}(A)$. When ${\rm ind}(A)\leq 1$, $A$ admits a group inverse, denoted by $A^\#$. Classical references on generalized inverses include \cite{benisraelgreville, campbell1979}, while more recent developments can be found in \cite{ChenBook}.  

Let $\psi_A$ and $\Delta_A$ denote the minimal and characteristic polynomials of $A$, respectively. It is well known that $\psi_A \mid \Delta_A \mid \psi^n_A$ and hence that $\psi_A$ and $\Delta_A$ share the same irreducible factors. The multiplicity of zero as a root of $\psi_A$ is called the \emph{index} of $A$, denoted $i(A)$; explicitly, if $\psi_A(\gl)=\gl^{i(A)}f(\gl)$ with $\gcd(\gl,f(\gl))=1$, then ${\rm ind}(A)=i(A)$.  
The Drazin index is invariant under similarity: since similar matrices share the same minimal polynomial, they necessarily have the same Drazin index. In particular, if $B=UAU^{-1}$ then 
\[
i(B)=i(A) \quad \text{and} \quad B^D = U A^D U^{-1}.
\]
Furthermore, for conformably dimensioned matrices $A,B$, one has $\psi_{AB}(\gl)=\gl^{0,\pm 1}\psi_{BA}(\gl)$ and hence $|i(AB)-i(BA)|\leq 1$. 

The Drazin inverse of a singular square matrix $A$ can be obtained through the so-called \emph{core--nilpotent decomposition}. Specifically, if $A$ is singular, then there exist invertible matrices $U$ and $C$, together with a nilpotent matrix $N$ satisfying $N^k=0 \neq N^{k-1}$ for some $k \geq 1$, such that 
\[
A = U \begin{bmatrix} C & 0 \\ 0 & N \end{bmatrix} U^{-1}.
\]
In this case,
\[
A^D = U \begin{bmatrix} C^{-1} & 0 \\ 0 & 0 \end{bmatrix} U^{-1},
\]
and ${\rm ind}(A)=k$.
The core--nilpotent decomposition immediately yields the following characterization:
\begin{lem} \label{lemagi}
Let $A$ be a square matrix. Then $i(A) = k$ if and only if $k$ is a positive integer such that $(A^k)^\#$ exists while
$(A^i)^\#$ do not, for all $0< i <k$.
\end{lem}
Furthermore, if $i(A)=k$ and $A^k = A^{k+1}X$ for some matrix $X$, then by applying the decomposition one obtains
\[
A^D = A^k X^{k+1}.
\]
This representation shows that the Drazin inverse can be expressed directly in terms of the minimal polynomial of $A$. Indeed, if $\psi_A(\gl)=\gl^k g(\gl)$ with $\gcd(\gl,g(\gl))=1$ and 
\[
g(\gl)=g_0 + g_1 \gl + \cdots + g_\ell \gl^\ell, \qquad \ell \geq 1,
\]
then from
\[
0 = A^k \big(g_0 I + g_1 A + \cdots + g_\ell A^\ell\big)
\]
it follows that
\[
g_0 A^k = A^{k+1}(-g_1 I - g_2 A - \cdots - g_\ell A^{\ell-1}),
\]
and hence
\[
A^k = A^{k+1}X, \quad \text{with} \quad X = -\tfrac{1}{g_0}(g_1 I + g_2 A + \cdots + g_\ell A^{\ell-1}).
\]
Cline, as shown in \cite{Cline}, established a fundamental relation between the Drazin inverses of $AB$ and $BA$, namely
\[
(AB)^D = A \big[(BA)^D\big]^2 B,
\]
which is commonly referred to as \emph{Cline's formula}.
The following theorem, also due to Cline \cite{Cline1968}, provides an alternative criterion for the existence of the group inverse together with an explicit formula for its computation.
\begin{thm}[Cline's Theorem]
Let $A$ be a square matrix with a full rank factorization $A = FG$. Then $A$ has a group inverse if and only if $GF$ is nonsingular. In that case,
\[
A^\# = F (GF)^{-2} G.
\]
\end{thm}
 A matrix $A$ is said to be Moore--Penrose invertible (with respect to $*$) if there exists a common solution $X$ to the following equations:
\begin{equation*}
\label{Eq:MoorePenrose}
\begin{aligned}
    &(1) \; A X A = A, \quad
    (2) \; XAX = X, \quad
    (3) \; (AX)^{*} = AX, \quad
    (4) \; (XA)^{*} = XA.
\end{aligned}
\end{equation*}
A solution to this set of equations is unique, if it exists, and is called the Moore--Penrose inverse of $A$, denoted by $A^\dag$.
Moore--Penrose invertibility can also be characterized via a full rank factorization, paralleling the group inverse case.
\begin{lem}[MacDuffee \cite{benisraelgreville}, Puystjens--Robinson \cite{Puystjens}] \label{macduffee} Let $A = FG$ be a full rank factorization of $A$. Then $A^\dagger$ exists if and only if both $GG^*$ and $F^*F$ are nonsingular. In this case,
\[
A^\dagger = G^*(GG^*)^{-1}(F^*F)^{-1}F^*.
\]
\end{lem}
Moore--Penrose invertibility is invariant under unitary equivalence. Indeed, if $A = U B V^*$ with $U^* = U^{-1}$ and $V^* = V^{-1}$, then $A^\dagger$ exists if and only if $B^\dagger$ exists. In that case,

\[
A^\dagger = V B^\dagger U^*.
\]

This property underlies many applications in numerical linear algebra and data science. 
A solution to equation (1), denoted by $A^{-}$, is referred to as a von Neumann inverse of $A$. The set of all such inverses is denoted by $A\{1\}$. A solution that simultaneously satisfies equations (1) and (2) is called a reflexive inverse of $A$, denoted by $A^{+}$. In particular, it is immediate that $A^{-} A A^{=}$ is a reflexive inverse of $A$ for any $A^{-}, A^{=} \in A\{1\}$. The set of common solutions to equations (1) and (3), which may be empty, is denoted by $A\{1,3\}$. An analogous notation is used for the common solutions to equations (1) and (4).

Let $G$ be a simple undirected graph on $n$ vertices $v_1,\ldots,v_n$. Its adjacency matrix is the $n\times n$ matrix $A=(a_{ij})$ given by
\[
a_{ij}=\begin{cases}
1 & \text{if } \{v_i,v_j\}\in E(G),\\
0 & \text{otherwise}.
\end{cases}
\]
The definition extends naturally to directed graphs (digraphs). Given a square matrix $A=[a_{ij}]$, its associated digraph $D(A)=(V,E)$ has vertex set $V=\{v_1,\ldots,v_n\}$ and edge set $E$, where $(v_i,v_j)\in E$ iff $a_{ij}\neq 0$. The entry $a_{ij}$ may be interpreted as the weight of the edge from $v_i$ to $v_j$, vanishing exactly when the edge is absent.  

Although the adjacency matrix depends on the ordering of the vertices, it is unique up to permutation similarity. If $A$ and $B$ are adjacency matrices of the same graph, then there exists a permutation matrix $P$ such that $A=P^{-1}BP$. Since $P^*=P^{-1}$, this relation constitutes a special case of unitary similarity. Consequently, the existence of Drazin, group, and Moore--Penrose inverses is independent of the chosen vertex ordering.  

In this paper we investigate Drazin and Moore--Penrose invertibility for matrices associated with two families of digraphs: in Section \ref{SecDoubleStar} we study double star digraphs, and in Section \ref{SecLinked} we analyze digraph linked stars, with particular emphasis on the role of the Drazin index and explicit inversion formulas.

\section{\textbf{Double star digraphs}}\label{SecDoubleStar}

\hspace{0.5cm} The {\it {double star digraph}} $S_{m,n}$ is obtained by connecting the central vertex of the star $K_{1,m-1}$ to the central vertex of the star $K_{1,n-1}$ with two directed edges, one in each direction.

Let $x$ and $y$ be strictly nonzero vectors of length $m$, while $z$ and $w$ are strictly nonzero vectors of length $n$. Consider the square matrix of order $m+n+2$ given by
\begin{equation}M=  \mxl{cc|cc}
0 & x^\top  & a & 0\\
y & 0 & 0 & 0\\ \hline
b & 0 & 0 & z^\top \\
0 & 0 & w & 0
\mxr,\label{matdsg}\end{equation}

where $a$ and $b$ are nonzero elements. Then $D(M)$ corresponds to a double star digraph $S_{(m+1),(n+1)}$ and any matrix $A$ with associated digraph $D(A)=S_{(m+1),(n+1)}$ is permutation similar to a matrix of the form described in (\ref{matdsg}).

In \cite{McDonald}, the authors present a formula for $M^\#$ in case $x^\top  y\neq 0$ and $z^\top  w\neq 0$. Note that $M^\#$ exists only in this case. In fact, $M$ admits the full rank factorization
\begin{equation}M= \mxl{cc|cc}
\label{fullrank}
0 & 1 & a & 0\\
y & 0 & 0 & 0\\ \hline
b & 0 & 0 & 1\\
0 & 0 & w & 0
\mxr \mxl{cc|cc}
1 & 0 & 0 & 0\\
0 & x^\top  & 0 & 0\\ \hline
0 & 0 & 1 &0\\
0 & 0 & 0 &  z^\top 
\mxr=FG,\end{equation} which gives $$GF= \bmx
0 & 1 & a & 0\\
x^\top  y & 0 & 0 & 0\\
b & 0 & 0 & 1\\
0 & 0 & z^\top  w & 0
\emx.$$
By Cline's Theorem, $M$ has a group inverse if and only if $GF$ is nonsingular, in which case $M^\#=F(GF)^{-2}G$. Since $\operatorname{rank}(GF)=4$ if and only if $x^\top y\neq 0$ and $z^\top w\neq 0$, $M^\#$ exists if and only if $x^\top y\neq 0$ and $z^\top w\neq 0$.

Having characterized the conditions for group invertibility, we now address the complementary situation in which this property fails.
In Subsection \ref{drazininverse}, we analyse the matrices $M$ associated with double star digraphs that are not group invertible, determining in each case the Drazin index and the minimal polynomial of $M$, and deriving an explicit expression for $M^D$. In Subsection \ref{moorepenrose} we investigate the Moore--Penrose invertibility of matrices over a field whose associated digraphs are of double star type, following the same structural framework. 
\; \\
\subsection{Drazin inverse}
\label{drazininverse}
\; \\
\; \\
\hspace{0.5cm} In this subsection we analyze all non group-invertible cases of matrices associated with double star digraphs. For each configuration determined by the relations among the vector and scalar entries of $M$, we identify the Drazin index, compute the minimal polynomial, and obtain a closed-form expression for the Drazin inverse $M^D$.

\begin{thm} \label{thm2.1} Let $M$ be a matrix over a field $\mathbb{F}$ with associated double star digraph $S_{(m+1),(n+1)}$ of the form (\ref{matdsg}), where $x^\top y=z^\top w=0$. Then $i(M)=2$, $\psi_{M}(\gl)=\gl^2 (\gl^2-ab)$, and
\begin{equation*} M^D= (ab)^{-1} \mxl{cc|cc}
0 & x^\top  & a & 0\\
y & 0 & 0 & b^{-1}yz^\top \\ \hline
b & 0 & 0 & z^\top \\
0 & a^{-1}wx^\top  & w & 0
\mxr,\label{matdsgthm1}\end{equation*}
\end{thm}

 \proof Consider the full rank factorization given in (\ref{fullrank}). Since $x^\top y=z^\top  w=0$, we can write
$$GF=\mxl{cc|cc}
0 & 1 & a & 0\\
0 & 0 & 0 & 0\\ \hline
b & 0 & 0 & 1\\
0 & 0 & 0 & 0\mxr.$$
Recall that $\psi_{FG}(\gl) = \gl^{0, \pm 1}\psi_{GF}(\gl)$ and $|i(FG)-i(GF)|\le 1$.
Let $M'=GF$. We show that $\psi_{M'}(\gl)=\psi_{GF}(\gl)=\gl (\gl^2-ab)$. Since $\lambda(\lambda-c)$ is not an annihilating polynomial of $M'$ for any $c$, it suffices to prove that $$ M'(M'^2-abI)=0 \quad \text{and} \quad M'^2 -abI\ne 0.$$ The equality follows directly by writing $$GF=\bmx X & Y\\ Z & W\emx$$ and observing that $X^{2}=W^{2}=XY=WZ=0$, $YZ=\operatorname{diag}(ab,0)$, $ZY=\operatorname{diag}(ba,0)$, $YW=\begin{bmatrix}0 & a\\0 & 0\end{bmatrix}$, and $ZX=\begin{bmatrix}0 & b\\0 & 0\end{bmatrix}$.
Hence $i(M')=1$, from which $i(M)\in \{0, 1, 2\}$. Since $i(M)=0$ would imply that $M$ is invertible, and $i(M)=1$ would imply that $M$ is group-invertible (which contradicts the singularity of $GF$), we conclude that
\[
i(M)=2.
\]
Therefore,
\[
\psi_{M}(\gl)=\psi_{FG}(\gl) = \gl\psi_{GF}(\gl)=\gl^2 (\gl^2-ab).
\]
From the minimal polynomial of $M$, it follows that $M^{2}=M^{3}X$ with 
\[
X = -\frac{1}{-ab}(0I+1M)=\frac{1}{ab}M,
\]
and consequently,
\begin{align*}
    M^D &= \frac{1}{(ab)^3}M^5 \\
    &= (ab)^{-1}\mxl{cc|cc}
0 & x^\top  & a & 0\\
y & 0 & 0 & b^{-1}yz^\top \\ \hline
b & 0 & 0 & z^\top \\
0 & a^{-1}wx^\top  & w & 0
\mxr.
\end{align*} \endproof

The next two results correspond to the case where $x^\top y\ne 0$ and $z^\top w=0$.

\begin{thm} \label{thm2.2} Let $M$ be a matrix over a field with associated double star digraph $S_{(m+1),(n+1)}$ of the form (\ref{matdsg}), where $x^\top y\ne 0$, $z^\top w=0$ and $\zeta = x^\top y+ab \ne 0$. Then $i(M)=3$, $\psi_{M}(\gl)=\gl^3 (\gl^2-\zeta)$, and
\begin{equation*}
    M^D= \zeta^{-1} \mxl{cc|cc}
0 & x^\top  & a & 0\\
y & 0 & 0 & \zeta^{-1}ayz^\top \\ \hline
b & 0 & 0 & \zeta^{-1}abz^\top \\
0 & \zeta^{-1}bwx^\top  & \zeta^{-1}abw & 0
\mxr.
\end{equation*}
\end{thm}

\proof Consider the full rank factorization given in (\ref{fullrank}).
Since $z^\top w=0$, we can write $$GF= \bmx
0 & 1 & a & 0\\
x^\top  y & 0 & 0 & 0\\
b & 0 & 0 & 1\\
0 & 0 & 0 & 0
\emx$$
Let $M'=GF$. Note that $M'$ admits the full rank factorization
$$M' =\bmx 0 & 1 & 0\\
x^\top  y & 0 &  0\\
b & 0  & 1\\
0 & 0 & 0
\emx \bmx 1 & 0 & 0 & 0\\ 0 & 1 & a & 0\\ 0 & 0 & 0 & 1\emx =F' G',$$ which gives
$$G'F' = \bmx 0 & 1 & 0\\\zeta & 0 & a\\ 0 & 0 & 0\emx.$$
Let $M''=G'F'$. Since
$$ M''(M''^2- \zeta I)=0 \ne M''^2 - \zeta I,$$
we conclude that $\psi_{M''}(\gl) =
 \gl(\gl^2-\zeta )$. Hence $\psi_{M'}(\gl)=\gl^j(\gl^2-\zeta )$ with $j\in\{0, 1, 2\}$. Because $M'$ is singular, $i(M') \neq 0$. By Lemma \ref{lemagi}, if $i(M') = 1$ then $M' = F'G'$ would have a group inverse and, by Cline's Theorem, $M'' = G'F'$ would be invertible -- contradicting the fact that $M''$ has a null row.  
Therefore, the cases $j = 0,1$ cannot occur, and $\psi_{M'}(\gl)=\gl^2(\gl^2-\zeta )$. This implies $\psi_M(\gl)=\gl^i (\gl^2-\zeta )$ with $i \in \{1,2,3\}$. 
Since $M'=GF$ is singular, then $M^\#$ cannot exist, and thus  $i(M)>1$. To show that $i(M)\ne 2$, note that
$$M^2 = \bmx \zeta & 0 &a\\
0 & y & 0\\
0 & b & 0\\
wb & 0 & w \emx
 \bmx 1 & 0 & 0 & 0\\
 0 & x^\top  & a & 0\\
 0 & 0 & 0 &z^\top  \emx = F_1 G_1$$ is a full rank factorization. Indeed, the first and third column of $F_1$ are linearly independent; otherwise, if $wb = w\mu$, then left-multiplying by $\tilde{w}$ gives $b = \mu$, and $\zeta = x^{T}y + ab = a\mu = ab$ would imply $x^{T}y = 0$, which is impossible.  
 If $i(M)=2$, then $(M^2)^\#$ would exist, implying that $G_{1}F_{1}$ is invertible.  Since $z^\top  w = 0$, we obtain $$G_1F_1 = \bmx \zeta & 0 & a\\
 0 & \zeta  & 0\\
 0 & 0 & 0\emx,$$ which is singular. Hence $i(M)\ne 2$.
 Therefore $i(M)=3$, and consequently $\psi_M(\gl)=\gl^3 (\gl^2-\zeta )$. Since $\psi_M(\gl)=\gl^3 (\gl^2-\zeta )$, we have $M^3=M^{4}X$, with
 \[
X = -\frac{1}{-\zeta}(0I + 1M) = \frac{1}{\zeta}M,
\]
and thus
\begin{align*}
M^D &= \frac{1}{\zeta^4}M^7\\
&= \zeta^{-1} \mxl{cc|cc}
0 & x^\top  & a & 0\\
y & 0 & 0 & \zeta^{-1}ayz^\top \\ \hline
b & 0 & 0 & \zeta^{-1}abz^\top \\
0 & \zeta^{-1}bwx^\top  & \zeta^{-1}abw & 0
\mxr.
\end{align*}
\endproof

\begin{thm} \label{thm2.3} Let $M$ be a matrix over a field with associated double star digraph $S_{(m+1),(n+1)}$ of the form (\ref{matdsg}), where $x^\top y=-ab$ and $z^\top w=0$. Then $i(M)=5$, $\psi_{M}(\gl)=\gl^5$, and $M^D=0$.
\end{thm}
 
 \proof Consider the full rank factorization given in (\ref{fullrank}).
Since  $x^\top y=-ab$ and $z^\top w=0$, we can write $$GF= \bmx
0 & 1 & a & 0\\
-ab & 0 & 0 & 0\\
b & 0 & 0 & 1\\
0 & 0 & 0 & 0
\emx.$$
Let $M'=GF$. Since $(M')^3\neq 0$ and $(M')^4=0$, we conclude that $\psi_{M'}(\gl)=\gl^4$ and, therefore, $\psi_{M}(\gl)=\gl^j$ with $j\in\{3,4,5\}$. 

Note that $M^4=(FG)^4=F(GF)^3G=F(M')^3G$. Hence
$$
M^4= \mxl{cc|cc}
0 & 0 & 0 & 0\\
0 & 0 & 0 & 0\\ \hline
0 & 0 & 0 &0\\
0 & 0 & 0 &  abwz^\top 
\mxr\neq 0.$$

This implies that $\psi_{M}(\gl)=\gl^5$ and $i(M)=5$. Therefore, $M$ is nilpotent of index 5, and consequently $M^D=0$.\endproof
Observe that the case where $x^\top y=0$ and $z^\top w\ne 0$ reduces, by permutation similarity, to the previous case where $x^\top y\ne 0$ and $z^\top w=0$.

The results presented in this subsection provide a complete description of the Drazin invertibility of matrices associated with double star digraphs in all non group-invertible cases. Having determined the Drazin index, minimal polynomial, and explicit expressions for $M^D$, we now turn our attention to Moore--Penrose invertibility within the same structural framework.

\subsection{Moore--Penrose inverse}
\label{moorepenrose}
\; \\
\; \\
\hspace{0.5cm} In this subsection, we investigate the Moore--Penrose invertibility of matrices associated with double star digraphs. We begin with a technical lemma that will serve as a key tool for deriving our main result.

\begin{lem} \label{lemaCondMP} The block matrix $\bmx \alpha & \beta \\ \gamma & 1 \emx$ is invertible if and only if  $\sigma = \alpha-\beta\gamma$ is invertible. In this case,
$$\bmx \alpha & \beta \\ \gamma & 1 \emx^{-1} = \bmx \sigma ^{-1} & -\sigma^{-1}\beta \\ -\gamma\sigma^{-1} & 1+ \gamma\sigma^{-1}\beta \emx.$$
\end{lem}
\proof The result follows directly by applying the Schur complement to the $(2,2)$ entry. \endproof
We are now in a position to establish necessary and sufficient conditions for the existence of the Moore--Penrose inverse of $M$, together with its explicit expression in terms of the component vectors.
\begin{thm} \label{thmMPdoublestar} Let $M$ be a matrix over a field $\mathbb{F}$ with associated double star digraph $S_{(m+1),(n+1)}$, of the form (\ref{matdsg}), where $x=[x_i], y = [y_i], z=[z_i]$ and $w=[w_i]$. Then
$M^\dagger$ exists if and only if $s,u,t,v\in \mathbb{F}\setminus \{0\}$, where 
\[
s = \sum x_i\overline{x_i}, \quad 
u = \sum y_i\overline{y_i}, \quad 
t = \sum z_i\overline{z_i}, \quad 
v = \sum w_i\overline{w_i}.
\]
In this case,
\begin{equation*}
    M^\dagger = \mxl{cc|cc} 0 & u^{-1}y^* & 0 & 0 \\ s^{-1}\bar{x} & 0 & 0 & -s^{-1}av^{-1}\bar{x}w^* \\ \hline 0 & 0 & 0 & v^{-1}w^* \\ 0 & -t^{-1}bu^{-1}\bar{z}y^* & t^{-1}\bar{z} & 0 \mxr.
\end{equation*}
\end{thm}
\proof
We start from the full rank factorization (\ref{fullrank}),
$$M=FG=\mxl{cc|cc}
0 & 1 & a & 0 \\
y & 0 & 0 & 0 \\
\hline
 b & 0 & 0 & \mathit{1} \\
0 & 0 & w & 0
\mxr \mxl{cc|cc}
1 & \mathit{0} & 0 & 0 \\
0 & \mathit{x^\top } & 0 & 0 \\
\hline
 0 & 0 & 1 & 0 \\
0 & 0 & 0 & \mathit{z^\top }
\mxr,$$   where $F$ is an $(m+n+2)\times 4$ matrix and $G$ is a $4\times (m+n+2)$ matrix, both of rank four.
Recalling Lemma \ref{macduffee}, $M^\dagger$ exists if and only if both $F^*F$ and $GG^*$ are nonsingular. We compute:
$$F^*F= \mxl{cc|cc}
0 & y^* & \overline{b} & 0 \\
1 & 0 & 0 & 0 \\
\hline
 \overline{a} & 0 & 0 & \mathit{w^*} \\
0 & 0 & 1 & 0
\mxr \mxl{cc|cc}
0 & 1 & a & 0 \\
y & 0 & 0 & 0 \\
\hline
 b & 0 & 0 & 1 \\
0 & 0 & w & 0
\mxr = \mxl{ccc|c}
y^*y+\overline{b}b & 0 & 0 & \overline{b} \\
0 & 1 & a & 0 \\
0 & \overline{a} & \overline{a}a+w^*w & 0 \\
 \hline
b & 0 & 0 & 1
\mxr.$$
Writing
$F^*F = \bmx \alpha & \beta \\ \gamma & 1 \emx$, with $\alpha = \mxl{ccc}
y^*y+\overline{b}b & 0 & 0\\
0 & 1 & a \\
0 & \overline{a} & \overline{a}a+w^*w \\
\mxr, \beta = \mx{\overline{b} \\ 0 \\ 0}$, and $\gamma = \bmx b & 0 & 0 \emx,$
and applying Lemma \ref{lemaCondMP}, $(F^*F)^{-1}$ exists if and only if $K=\alpha-\beta \gamma$ is invertible. In that case,
\begin{equation*} \label{Eq.F*F}
    (F^*F)^{-1} = \bmx K^{-1} & -K^{-1}\beta \\ -\gamma K^{-1} & 1+\gamma K^{-1}\beta \emx.
\end{equation*}
Note that
$K =\mxl{ccc}
y^*y & 0 & 0\\
0 & 1 & a \\
0 & \overline{a} & \overline{a}a+w^*w \\
\mxr$.
We may take the factorization

$PKP=\mxl{cc|c} y^*y & 0 & 0\\ 0 & \overline{a}a+w^*w & \overline{a} \\ \hline 0 & a & 1 \\ \mxr=X$, with   $P = \mxl{ccc} 1 & 0 & 0 \\ 0 & 0 & 1 \\ 0 & 1 & 0 \mxr$. Obviously $K$ is invertible if and only if $X$ is invertible, in which case  $K^{-1}=PX^{-1}P$.

Applying Lemma \ref{lemaCondMP} again, $X$ is invertible if and only if $\mxl{cc} y^*y & 0\\ 0 & \overline{a}a+w^*w \\ \mxr - \mxl{c}0 \\ \overline{a}\mxr \mxl{cc} 0 & a\mxr$ is invertible. That is, $X$ is invertible exactly when  $$S=\mxl{cc} y^*y & 0\\ 0 & w^*w \\ \mxr $$ is   invertible. This occurs precisely when $u=y^*y$ and $v=w^*w$ are nonzero elements of $\mathbb{F}.$
Now $$S^{-1}=\mxl{cc} u^{-1} & 0\\ 0 & v^{-1} \\ \mxr \quad
\text{and by Lemma \ref{lemaCondMP},}\quad X^{-1}=\mxl{ccc} u^{-1} & 0 & 0\\ 0 & v^{-1} & -v^{-1}\overline{a} \\ 0 & -av^{-1} & 1+av^{-1}\bar{a} \mxr,$$ from which
$$
K^{-1} = P X^{-1} P = \mxl{ccc} u^{-1} & 0 & 0 \\ 0 & 1+ av^{-1}\overline{a} & -a v^{-1} \\ 0 & -\overline{a}v^{-1} & v^{-1} \mxr.$$
Applying Lemma \ref{lemaCondMP} once more, we obtain
$$(F^*F)^{-1}=\bmx K^{-1} & -K^{-1}\beta \\ -\gamma K^{-1} & 1+\gamma K^{-1}\beta \emx =  \mxl{cccc} u^{-1} & 0 & 0 & -\bar{b} u^{-1} \\ 0 & 1+ a \overline{a}v^{-1} & -a v^{-1} & 0 \\ 0 & -\overline{a}v^{-1} & v^{-1} & 0 \\ -bu^{-1} & 0 & 0 & 1+b \bar{b}u^{-1} \mxr.$$
Considering now $GG^*$, we have
$$GG^* = \mxl{cc|cc} 1 & 0 & 0 & 0 \\ 0 & \mathit{x^\top } & 0 & 0 \\ \hline 0 & 0 & 1 & 0 \\ 0 & 0 & 0 & \mathinner{z^\top } \mxr \mxl{cc|cc} 1 & 0 & 0 & 0 \\ 0 & \mathit{\bar{x}} & 0 & 0 \\ \hline 0 & 0 & 1 & 0 \\ 0 & 0 & 0 & \mathinner{\bar{z}} \mxr = \mxl{cc|cc} 1 & 0 & 0 & 0 \\ 0 & \mathit{x^\top  \overline{x}} & 0 & 0 \\ \hline 0 & 0 & 1 & 0 \\ 0 & 0 & 0 & \mathinner{z^\top  \overline{z}} \mxr.$$
Therefore, $GG^*$ is invertible if and only if $s=x^\top \overline{x} \neq 0$ and  $t=z^\top \overline{z} \neq 0$. In this case, $$(GG^*)^{-1} = \mxl{cc|cc} 1 & 0 & 0 & 0 \\ 0 & s^{-1} & 0 & 0 \\ \hline 0 & 0 & 1 & 0 \\ 0 & 0 & 0 & t^{-1} \mxr.$$
Finally, by Lemma~\ref{macduffee}, the expression for $M^\dagger$ follows as stated.\endproof
\section{\textbf{$D$-linked stars digraphs}}\label{SecLinked}
\hspace{0.5cm} In this section we extend our previous analysis to a broader class of digraphs, namely the
\emph{D-linked stars digraphs}. This family generalizes the double star digraphs discussed in Section \ref{SecDoubleStar},
allowing for an arbitrary linking structure determined by a base digraph $D$. As before, our goal is to determine explicit conditions for the existence of the Drazin and Moore--Penrose inverses of the corresponding matrices, together with closed-form expressions.

We begin by introducing the construction and recalling the notation to be used throughout this section.

Let $D$ be a digraph on vertices $v_1, v_2, \ldots, v_n$, and consider the directed stars $K_{1,r_1}, K_{1,r_2},$ $\ldots,K_{1,r_n}$. Following the definition in \cite{McDonald}, the \emph{$D$-linked stars digraph}, denoted $\operatorname{gls}(D,r_1,r_2,$ $\ldots,r_n)$, is obtained by connecting the centre of $K_{1,r_i}$ to the centre of $K_{1,r_j}$ whenever $(v_i,v_j)$ is an edge of $D$.
Note that double star digraphs constitute a particular case of D-linked stars.

For an $n\xx n$ matrix $A$ over $\FF$ and strictly nonzero vectors $x_i, y_i$ of length $r_i\in\mathbb{N}$,
consider the block matrix
\begin{equation}
M =
\begin{bmatrix}
A & B \\[2mm]
C & 0
\end{bmatrix},
\qquad
B = \operatorname{diag}(x_1^T,\ldots,x_n^T), \quad
C = \operatorname{diag}(y_1,\ldots,y_n).
\label{matgls}
\end{equation}

Then $D(M)$ is a digraph of type D-linked stars associated with $D(A)$. Any matrix associated with such a digraph is permutation similar to one of the form (\ref{matgls}). Observe that the matrix corresponding to $D$-linked stars is in particular a (220) matrix, as defined in \cite{Grp220}.

We next examine, separately, the conditions for Drazin and Moore--Penrose invertibility.
\; \\
\subsection{Drazin inverse}
\label{subsecDlinked_Drazin}
\; \\
\; \\
\hspace{0.5cm} We first address the Drazin (and group) invertibility of matrices associated with $D$-linked stars digraphs.
Observe that the matrix $M$ in (\ref{matgls}) admits the full rank factorization \begin{equation}M=\bmx A & B\\ C & 0\emx = \bmx A & I\\ C & 0\emx \bmx I & 0\\ 0 & B\emx=FG,\label{fullrankgls}\end{equation} which yields $$GF= \bmx A & I\\ BC & 0\emx.$$

By Cline's Theorem, $M$ has a group inverse if and only if $GF$ is nonsingular. We now formalize this criterion.

\begin{thm}\label{4.1} Let $A=[a_{ij}]$ be an $n\times n$ matrix over an arbitrary field, and let $x_i, y_i$ be strictly nonzero vectors of length $r_i$. For $M$ defined as in (3), the following conditions are equivalent
\begin{enumerate}
\item $M$ has a group inverse;
\item $BC$ is invertible;
\item $x_i^\top y_i\ne 0$ for all $i=1,\dots,n$.
\end{enumerate}
Moreover,
$$M^\#=\bmx
0 & (BC)^{-1}B\\
C(BC)^{-1} & -C(BC)^{-1}A(BC)^{-1}B\emx.$$
\end{thm}

\proof
By the factorization (\ref{fullrankgls}), $M$ has a group inverse if and only if $GF$ is invertible,
which occurs precisely when $BC$ is invertible.
Since $BC = \operatorname{diag}(x_i^\top  y_i)$, the stated equivalences follow directly.
The explicit formula for $M^{\#}$ follows from Cline's theorem.
\endproof

When $BC$ is singular, the group inverse does not exist; it is therefore natural to investigate the Drazin invertibility in this case.

For matrices \(B\) and \(C\) as in (\ref{matgls}), \(B\) is of full row rank and \(C\) is of full column rank, both equal to $n$. Hence, there exist a right inverse of \(B\) and a left inverse of \(C\), denoted by \(B^{+}\) and \(C^{+}\), respectively. These matrices serve as reflexive inverses of \(B\) and \(C\); in particular, \(BB^{+} = I\) and \(C^{+}C = I\), implying that \(B^{+} \in B\{1,3\}\) and \(C^{+} \in C\{1,4\}\). From this point on, the notation $B^+$ and $C^+$ will consistently refer to these right and left inverses of $B$ and $C$, respectively.

In the following result, we consider the case where $BC=0$ and we determine $i(M)$ in terms of $i(A)$. 
Knowing $\psi_{A}(\gl)$ allows us to determine $i(M)$ and $\psi_{M}(\gl)$, and consequently to express $M^D$ in terms of the minimal polynomial.

\begin{thm}\label{thmDlinked_Drazin}
Let $A=[a_{ij}]$ be an $n\times n$ matrix over $\mathbb{F}$,
and let $x_i, y_i$ be strictly nonzero vectors of length~$r_i$ such that $x_i^\top  y_i = 0$ for all~$i$.
If $M$ is defined as in (\ref{matgls}), then
\[
i(M) = i(A) + 2.
\]
\end{thm}

\proof
Consider the full rank factorization (\ref{fullrankgls}).  
Since $x_i^\top  y_i = 0$ for all $i$, we have $BC=0$, so that
\[
GF = \begin{bmatrix} A & I \\[1mm] 0 & 0 \end{bmatrix}.
\]
Denote $M' = GF$.  By induction one verifies that, for every integer $j\ge1$,
$$
(M')^j = \begin{bmatrix} A^j & A^{j-1} \\[1mm] 0 & 0 \end{bmatrix},$$
and, for every integer $j\ge2$,
$$M^j = \begin{bmatrix} A^j & A^{j-1}B \\[1mm] CA^{j-1} & CA^{j-2}B \end{bmatrix}.$$
Suppose that $\psi_A(\lambda) = \lambda^{k} g(\lambda)$, where $\lambda \nmid g(\lambda)$ and 
\[
g(\lambda) = g_0 + g_1\lambda + \cdots + g_\ell \lambda^\ell, \qquad (\ell \ge 1).
\]
Therefore,
\[
A^{k}(g_0 I + g_1 A + \cdots + g_\ell A^\ell) = 0,
\]
and $k$ is the smallest nonnegative integer for which this equality holds.
We now prove that $\psi_{M'}(\lambda) = \lambda^{k+1} g(\lambda)$.  
To this end, we first show that $\lambda^{k+1} g(\lambda)$ annihilates $M'$.
Indeed,
\begin{eqnarray*}
&& M'^{\,k+1}(g_0 I + g_1 M' + \cdots + g_\ell M'^\ell)
\\[1ex]
&=& M'^{\,k+1}
\left(
g_0 \mxl{c|c} I & 0 \\ \hline 0 & I \mxr
+ g_1 \mxl{c|c} A & I \\ \hline 0 & 0 \mxr
+ \cdots
+ g_\ell \mxl{c|c} A^\ell & A^{\ell-1} \\ \hline 0 & 0 \mxr
\right)
\\[1ex]
&=&
\mxl{c|c}
A^{k+1}(g_0 I + g_1 A + \cdots + g_\ell A^\ell) &
A^{k}(g_0 I + g_1 A + g_2 A^2 + \cdots + g_\ell A^{\ell}) \\ \hline
0 & 0
\mxr,
\end{eqnarray*}
which equals the zero matrix.
Hence, $\lambda^{k+1} g(\lambda)$ annihilates $M'$.
To show that $\lambda^{k} g(\lambda)$ does \emph{not} annihilate $M'$, 
note that if it did, we would have
\[
A^{k-1}(g_0 I + g_1 A + \cdots + g_\ell A^\ell) = 0,
\]
which contradicts the minimality of $k$.  
Therefore, $\psi_{M'}(\lambda) = \lambda^{k+1} g(\lambda)$.
From the above, we conclude that
\[
i(M') = k + 1 = i(A) + 1,
\]
and consequently $i(M) \in \{k,\, k+1,\, k+2\}$.
Assume, for contradiction, that $i(M) = k + 1$.  
Then
\begin{align*}
0 
&= M^{k+1} g(M) \\[1ex]
&= 
\bmx
A^{k+1} & A^k B \\[1mm]
C A^k & C A^{k-1} B
\emx\times\\
&\times\bmx
g_0 I + g_1 A + \cdots + g_\ell A^\ell &
g_1 B + g_2 A B + g_3 A^2 B + \cdots + g_\ell A^{\ell-1} B \\[1mm]
g_1 C + g_2 C A + \cdots + g_\ell C A^{\ell-1} &
g_0 I + g_2 C B + g_3 C A B + \cdots + g_\ell C A^{\ell-2} B
\emx.
\end{align*}
Since $BC = 0$, the $(2,2)$ block of this product gives
\[
C\!\left(A^{k-1}(g_0 I + g_1 A + \cdots + g_\ell A^\ell)\right) B = 0.
\]
Multiplying on the left by $C^{+}$ and on the right by $B^{+}$ yields
\[
A^{k-1}(g_0 I + g_1 A + \cdots + g_\ell A^\ell) = 0,
\]
which means that $\lambda^{k-1} g(\lambda)$ is an annihilating polynomial of $A$.
This contradicts the minimality of $\psi_A(\lambda)$.
Hence, the assumption $i(M) = k + 1$ is false, and we conclude that
$i(M) = i(A) + 2$. $\quad$ \endproof

The results in Theorems \ref{4.1} and \ref{thmDlinked_Drazin}
characterize group invertibility and determine the Drazin index for matrices associated with $D$-linked stars digraphs. We now turn to the Moore--Penrose inverse.

\; \\
\subsection{Moore--Penrose inverse}\label{subsecDlinked_MP}
\; \\
\; \\
\hspace{0.5cm} We next determine conditions under which the matrix $M$ in (\ref{matgls})
admits a Moore--Penrose inverse and derive its explicit form, relying on a suitable full rank factorization.

\begin{thm}
Let $A=[a_{ij}]$ be an $n\times n$
 matrix over a field $\FF$. Let $x_i, y_i$ be strictly nonzero vectors of length $r_i$, and let $M$ be defined as in (\ref{matgls}). Then $M^\dagger$ exists if and only if $x_i^*x_i, y_i y_i^*\in \FF\setminus\{0\}$, for all $i$. In this case, $$M^\dagger = \bmx 0 & (C^*C)^{-1}C^*  \\  B^*(BB^*)^{-1} &  -B^*(BB^*)^{-1} A (C^*C)^{-1}C^* \emx.$$
\end{thm}

\proof
The block matrix $M = \bmx A & B \\ C & 0 \emx$ is unitarily equivalent to
\[
L = \bmx C & 0 \\ A & B \emx,
\]
which admits the full rank factorization $L = FG$, where
\[
F = \bmx C & 0 \\ A & I \emx,
\qquad
G = \bmx I & 0 \\ 0 & B \emx.
\]
According to Lemma \ref{macduffee}, $L$ (and hence $M$) is Moore--Penrose invertible
if and only if both $F^*F$ and $GG^*$ are nonsingular. This is equivalent to the invertibility of $C^*C$ and $BB^*$, which occurs precisely when $x_i^*x_i,\, y_i y_i^* \neq 0$ for all $i$.
Since $B^*(BB^*)^{-1}\in B\{1, 4\}$ and $B^+\in B\{1, 3\}$, it follows that $$B^\dagger = B^*(BB^*)^{-1}BB^+ = B^*(BB^*)^{-1}.$$ Analogously, $$C^\dagger = (C^*C)^{-1}C^*.$$
The equality $A = BB^+ AC^+ C $ implies that $BB^\dagger A = A$, and hence $BB^+ AC^+ = AC^\dagger$.  Similarly, we obtain $AC^\dagger C=A$. From these identities, it follows that $$L^\dagger = \bmx C^\dagger & 0 \\ -B^\dagger A C^\dagger & B^\dagger \emx,$$ from which
  $$M^\dagger =\bmx 0 & C^\dagger  \\  B^\dagger &  -B^\dagger A C^\dagger \emx =\bmx 0 & (C^*C)^{-1}C^*  \\  B^*(BB^*)^{-1} &  -B^*(BB^*)^{-1} A (C^*C)^{-1}C^* \emx.$$
\endproof

\section*{\textbf{Questions and remarks}}

\hspace{0.5cm} The results presented in this paper provide explicit and unified formulations for the Drazin and Moore--Penrose inverses of matrices associated with double star and $D$-linked star digraphs. The analysis clarifies how the Drazin index and minimal polynomial reflect the combinatorial structure of these graphs, bridging algebraic and spectral properties within a coherent framework. Several open questions and research directions naturally emerge from our findings:

\begin{enumerate}
\item
The existence and explicit form of the group and Moore--Penrose inverses of matrices associated with $D$-linked stars digraphs 
have been completely characterized. However, in the case of the Drazin inverse, we only determined the Drazin index when $BC=0$,
and an explicit closed-form expression for $M^D$ remains to be found, even in this particular case.
A natural question is whether $M^D$ can be expressed constructively
in terms of $A^D$ and the off-diagonal blocks $B$ and $C$.

\item
The results obtained so far concern matrices over fields endowed with an involution. Extending this framework to rings or *-algebras could reveal deeper algebraic properties and broaden the applicability of generalized inverses in noncommutative settings.
\item
Theorem \ref{thmDlinked_Drazin} relates the Drazin indices of $A$ and $M$ when $BC=0$.
It remains to determine whether a full spectral correspondence
\[
\sigma(M)\setminus\{0\} = \sigma(A)\setminus\{0\}
\]
holds under general conditions, or whether additional structural assumptions on $B$ and $C$ are required.
\item The block formulas derived in Sections \ref{moorepenrose} and \ref{subsecDlinked_MP} are suitable for symbolic or algorithmic computation. A detailed study of their numerical stability, complexity, and performance for large sparse matrices would be valuable, particularly in connection with iterative or perturbation-based methods for computing generalized inverses.
\item Extending the present approach to other structured digraphs -- such as block-tree, multipartite, or hierarchically linked networks -- may yield new insights into how algebraic invariants encode graph topology.
\item Finally, one may ask whether, given a prescribed Drazin index or minimal polynomial, it is possible to construct a digraph (or a family of adjacency blocks) realizing such properties. Addressing this ``inverse spectral problem'' could enrich the combinatorial interpretation of generalized inverses.
\end{enumerate}

\section*{Acknowledgements}
{This research was partially financed by Portuguese Funds
through FCT (Funda\c c\~ao para a Ci\^encia e a Tecnologia) within the Project UID/00013/2025.}
\section*{Data Availability}
No datasets were generated or analyzed during the current study. All results are theoretical and fully contained within the manuscript.
\section*{Conflict of interest}
The authors declare that they have no conflict of interest.

\bibliographystyle{plain}
\bibliography{Bibliography}
\end{document}